\documentclass[12pt]{article}
\newcommand{\al}{\alpha}
\newcommand{\be}{\beta}
\newcommand{\ga}{\gamma}
\newcommand{\pddt}[1]{\frac{\partial #1}{\partial t}}

\begin{document}
\title{Singularity formation in the Yang-Mills Flow}
\newcounter{theor}
\setcounter{theor}{1}
\newtheorem{theorem}{Theorem}[section]
\newtheorem{lemma}{Lemma}[section]
\newtheorem{remark}{Remark}[section]
\newtheorem{corollary}{Corollary}[section]
\newenvironment{proof}[1][Proof]{\begin{trivlist}
\item[\hskip \labelsep {\it #1}]}{\end{trivlist}}

\centerline{\bf SINGULARITY FORMATION IN THE YANG-MILLS FLOW}
\bigskip
\centerline{Ben Weinkove}
\centerline{Department of Mathematics, Columbia University}
\centerline{New York, NY 10027}
\centerline{{\it Email:} weinkove@math.columbia.edu}

\bigskip
\noindent
{\it Mathematics Subject Classification: 53C44}

\setcounter{section}{1}
\bigskip
\bigskip
\noindent
{\bf 1.  Introduction}
\bigskip

The Yang-Mills heat flow was shown to exist for all time over 
compact K\"ahler 
surfaces by Donaldson \cite{do} (see also \cite{dk}).  Over a compact 
Riemannian manifold, R{\aa}de \cite{ra} proved that for
dimensions two and three,  the flow exists for all time and
converges to a Yang-Mills connection.   However, his methods
fail for higher dimensions.  In dimensions $n>4$ it has  been
shown by Naito
\cite{na} that singularities can occur in finite time  (see
also \cite{gr}).  The  case of $n=4$ is still unknown (see
\cite{sst} and \cite{st}). The Yang-Mills heat flow and the
harmonic map heat flow appear to share similar  properties.  In
\cite{gh}, Grayson and Hamilton show using a monotonicity 
formula from \cite{ha} that rapidly  forming singularities in
the harmonic map heat flow converge, after a  blow-up process,
to homothetically shrinking solitons.  

The same technique works in the case of the Yang-Mills flow.  There is, 
however, a complication arising from the presence of the gauge group.  For this reason, we
must allow the freedom to make gauge transformations.   We use Hamilton's monotonicity
formula for the Yang-Mills flow
\cite{ha} to  prove that a sequence of blow-ups of a rapidly forming singularity will 
converge, modulo the gauge group, to a non-trivial  
homothetically shrinking soliton.  In the last section, explicit examples 
of such solitons are given in the case of trivial bundles over $R^n$ for 
$5 \le n \le 9$.

\bigbreak

We will first introduce some notation.
Let $M$ be a compact $n$-dimensional Riemannian manifold.  Let $E$ be a 
real
vector 
bundle of rank $r$ over $M$ with compact structure group $G$ and assume the 
fibres carry an inner product and that $G \subset SO(r)$ respects this inner 
product. 

We will often work in local coordinates, using Greek letters for 
the bundle indices and Latin letters for the 
manifold indices.  A connection $A$ is locally an 
endomorphism valued 1-form given by 
$A^{\al}_{i \be}$.  The curvature of the connection is given by
$$F^{\al}_{ij \be} = \partial_i A^{\al}_{j \be} - \partial_j
A^{\al}_{i \be} + A^{\al}_{i \gamma} A^{\gamma}_{j \be} - A^{\al}_{j
  \gamma} A^{\gamma}_{i \be}.$$

The Yang-Mills heat flow is a flow of connections $A=A(t)$ given by
$$ \frac{\partial}{\partial t} A^{\al}_{j \be} = D_p F^{\al}_{p j \be},$$
where $D$ denotes covariant differentiation with respect to the connection.  In
normal coordinates, this right hand side is given by
$$D_p F^{\al}_{p j \be} = \partial_p F^{\al}_{p j \be} + A^{\al}_{p \ga} 
F^{\ga}_{p j \be}
- F^{\al}_{p j \ga} A^{\ga}_{p \be}.$$ 
The Yang-Mills flow is a gradient flow for the Yang-Mills
functional
$$E = \frac{1}{2} \int_M |F|^2.$$ 
Let $s = s^{\al}_{\ \be}$ be a gauge transformation.  Then $s$ acts on 
a connection matrix $A$  by
$$s \cdot A = sAs^{-1} - (ds)s^{-1}.$$

Define a flow of connections $A=A(x,t)$ on the trivial bundle over 
$\mathbf{R}^n \times (-\infty,0)$ to be a homothetically shrinking soliton if 
it is a solution to the Yang-Mills heat equation and
satisfies the dilation condition
$$A^{\al}_{i \be} (x,t) = \lambda A^{\al}_{i \be} (\lambda x, \lambda^2 
t).$$

The main theorem is as follows.

\bigskip
\noindent
{\bf Main Theorem.} {\it Let $A$ be a smooth solution of the Yang-Mills 
flow on $M$ for $0 
\le t <T$ with a singularity at some point $X$ as $t \rightarrow T$.  Suppose 
that it is a rapidly forming singularity, so that 
$$ (T-t) |F| \le C,$$ for some constant $C$.  Then a sequence
of blow-ups around $(X,T)$,
$$A_p (i) (x,t) = \lambda_i A_p (\lambda_i x , T+ \lambda_i^2 t),$$
with factor $\lambda_i \rightarrow 0$, has a
subsequence which converges in
$C^{\infty}$ on compact sets modulo the gauge group to a 
homothetically shrinking soliton
 $\overline{A}$ defined on $\mathbf{R}^n \times (-\infty,0)$ which has 
non-zero curvature.}
\bigskip

In the course of the proof of the theorem, the meaning of the phrase 'modulo the gauge
group' will be made more precise.

\bigskip

The outline of the paper is as follows. In section 2, we show that given 
a bound on the curvature, we can derive bounds on all of the derivatives 
of the curvature.  In section 3, we use these 
estimates, together with a theorem of Uhlenbeck's \cite{u82}, to get 
bounds on the connections and then
convergence of a sequence of blow-ups near the singular point.  We  
use Hamilton's monotonicity formula \cite{ha} to show that, allowing gauge transformations,
the  sequence converges to a homothetically shrinking soliton.  We then use
an `$\epsilon$-regularity' result to show that this soliton must have 
non-zero curvature.  Finally, in section 4, we give examples of 
homothetically shrinking solitons on a trival $SO(n)$ bundle over 
$\mathbf{R}^n$ for $5 \le n \le 9.$

The results of this paper will form part of my forthcoming PhD thesis at
Columbia University. I would like to thank my thesis advisor D.H. Phong,
for his constant support and advice, and
Richard Hamilton for suggesting an improvement to an initial draft of the result.
I am also grateful to the referee for pointing out a couple of errors in an earlier version
of this paper.

\pagebreak[3]

\bigskip
\bigskip
\noindent
\setcounter{section}{2}
{\bf 2. Estimates on derivatives of the curvature}
\bigskip

We begin by proving local estimates on the derivatives of the curvature.  
First define the parabolic cylinder
$$P_r (X, T) = \{ (x,t) \in M \times \mathbf{R} \ | \ d(X,x) \le r, \ T-r^2 \le 
t \le T \}.$$
The following lemma is proved in \cite{gh} (the conclusion is slightly 
different, but it follows easily.)

\pagebreak[3]
\addtocounter{theorem}{1}
\begin{lemma} \label{lemmah}
There exists a constant $s>0$ and for every $\ga <1$, a constant $C_{\ga}$, 
such that if $h$ is a smooth function on $M$ with
$$ \pddt{h} \le \triangle h - h^2$$
when $h \ge 0$ in some parabolic cylinder $P_r (X,T)$ with $r\le s$ then
$$h \le \frac{C_{\ga}}{r^2}$$
on $P_{\ga r} (X,T)$.
\end{lemma}

Then we have 

\addtocounter{lemma}{1}
\begin{theorem} \label{theoremestimates}
There exists a constant $s>0$ and constants $C_k$ for $k=1,2,\ldots$ depending 
only on $M$ such that if $A(t)$ is a solution to the Yang-Mills flow with $|F| 
\le K$ in some $P_r(X,T)$ for $r\le s$ for some constant $K\ge 1/r^2$, then 
with $r_k = r/2^{k}$ we have
$$ |D^k F| \le C_k K^{k/2+1} \qquad \textrm{on} \ P_{r_k}(X,T).$$
\end{theorem}
\begin{proof}
We will use induction to prove the theorem.  The notation $S * T$ 
will 
mean some algebraic bilinear expression involving the tensors $S$ and $T$. 
$R$ will denote the Riemannian curvature tensor on the base manifold. $C$ 
will denote a constant which depends only on $M$.  Assume that the 
constant $K$ is at least 1.  We begin with the case $k=1$.  Using 
commutation formulae and the Bianchi identity, calculate
\begin{eqnarray*}
\pddt{} |F|^2 & = & \triangle |F|^2 - 2|DF|^2 + F * F * F + R * F * F \\
& \le & \triangle |F|^2 - 2|DF|^2 + CK^3,
\end{eqnarray*}
and
\begin{eqnarray*}
\pddt{} |DF|^2 & = & \triangle |DF|^2 - 2|D^2 F|^2 + DF * DF * F \\
&& \mbox{} + DF * DF*R + DF* F * DR \\
& \le & \triangle |DF|^2 - 2|D^2 F|^2 + CK |DF|^2 + CK|DF| \\
& \le & \triangle |DF|^2 - 2|D^2 F|^2 + CK|DF|^2 + CK^4.
\end{eqnarray*}

Define
$$h = (8K^2 + |F|^2) |DF|^2.$$
Then 
\begin{eqnarray*}
\pddt{h} & \le &(\triangle|F|^2 - 2|DF|^2 + CK^3) |DF|^2 + (8K^2 + |F|^2) 
(\triangle |DF|^2 - 2|D^2 F|^2 \\
&& \mbox{} + CK|DF|^2 +CK^4).
\end{eqnarray*}
But
\begin{eqnarray*}
\triangle h \ge (\triangle |F|^2 ) |DF|^2 + (8K^2 + |F|^2) \triangle |DF|^2 - 
8K |DF|^2 |D^2 F|.
\end{eqnarray*}
Hence
\begin{eqnarray*}
\pddt{h} & \le & \triangle h - 2|DF|^4 + CK^3 |DF|^2 + (8K^2 + |F|^2) (-2|D^2 
F|^2 \\
&& \mbox{} + CK |DF|^2 + CK^4) + 8K |DF|^2 |D^2 F| \\
& \le & \triangle h - 2|DF|^4 + CK^3 |DF|^2 + CK^6  - 2(8K^2 + |F|^2) |D^2 
F|^2\\
& & \mbox{} + |DF|^4 + 16K^2 |D^2 F|^2 \\
& \le & \triangle h - \frac{1}{2} |DF|^4 + CK^6 \\
& \le & \triangle h - \frac{h^2}{C K^4} + CK^6
\end{eqnarray*}
for some constant $C$.  For that same $C$, define
$$ \hat{h} = \frac{h}{CK^4} - K.$$
Then for $\hat{h} \ge 0$,
\begin{eqnarray*}
\pddt{\hat{h}} & \le & \frac{1}{CK^4} (\triangle h - \frac{h^2}{CK^4} + 
CK^6) \\
& = & \triangle \hat{h} - (\hat{h} + K)^2 + K^2 \\
& \le & \triangle \hat{h} - \hat{h}^2.
\end{eqnarray*}
Then by Lemma \ref{lemmah} we have
$$ \hat{h} \le \frac{C}{r^2} \le C K$$
on $P_{r/2} (X,T)$.  Hence
$$h \le CK^5$$
and 
$$|DF| \le CK^{3/2}.$$
Assume inductively that we have the estimates for $ D^l F$ on 
$P_{r_{k-1}}(X,T)$ for $1 \le l \le k-1$.  Calculating on 
$P_{r_{k-1}}(X,T)$, first notice
\begin{eqnarray*}
\pddt{} |D^2 F|^2 & = & \triangle |D^2 F|^2 - 2|D^3 F|^2 + D^2 F * D^2 F * 
F + D^2 F * DF * DF \\
&& \mbox{} + D^2 F * D^2F * R + D^2 F * DF * DR + D^2 F * F * D^2 R \\
& \le & \triangle |D^2 F|^2 - 2|D^3 F|^2 + CK^5.
\end{eqnarray*}
It is not difficult to see by an induction argument that
\begin{eqnarray*}
\pddt{} |D^{k-1} F|^2 \le \triangle |D^{k-1} F|^2 - 2|D^{k}F|^2 + 
CK^{k+2}.
\end{eqnarray*}
and
\begin{eqnarray*}
\pddt{} |D^{k} F|^2 \le \triangle |D^{k} F|^2 - 2|D^{k+1}F|^2 + CK|D^kF|^2 
+ CK^{k+3}.
\end{eqnarray*}
Choose $B$ to be a constant with $K^{(k+1)/2} \le B \le CK^{(k+1)/2} $ and 
$|D^{k-1}F| \le B$.  Define
$h = (8 B^2 + |D^{k-1}F|^2) |D^k F|^2$.
Then 
\begin{eqnarray*}
\pddt{h} & \le & (8B^2 + |D^{k-1}F|^2) (\triangle |D^k F|^2 - 2|D^{k+1} 
F|^2 + CK |D^k F|^2 + CK^{k+3}) \\
&& \mbox{} + |D^k F|^2 (\triangle |D^{k-1}F|^2 - 2 |D^k F|^2 + CK^{k+2}) \\ 
& \le & (8B^2 + |D^{k-1}F|^2) \triangle |D^k F|^2 + |D^k F|^2 \triangle 
|D^{k-1}F|^2 - 16B^2 |D^{k+1} F|^2 \\
&& \mbox{} - 2 |D^k F|^4  + C K^{k+2} |D^k F|^2 + CK^{2k+4}.
\end{eqnarray*}
But
$$ \triangle h \ge (8B^2 + |D^{k-1}F|^2) \triangle |D^k F|^2 + |D^k F|^2 
\triangle |D^{k-1}F|^2 - 8B |D^k F|^2 |D^{k+1} F|.$$
Hence
\begin{eqnarray*}
\pddt{h} & \le & \triangle h - 16 B^2 |D^{k+1} F|^2 - 2|D^k F|^4 + 
CK^{k+2} |D^k F|^2 + CK^{2k+4}\\
&& \mbox{} + 16B^2 |D^{k+1} F|^2 + |D^k F|^4\\
& \le & \triangle h - |D^k F|^4 + CK^{k+2} |D^k F|^2 + CK^{2k+4} \\
& \le & \triangle h - \frac{1}{2} |D^k F|^4 + CK^{2k+4} \\
& \le & \triangle h - \frac{h^2}{CK^{2k+2}} + CK^{2k+4}.
\end{eqnarray*}
Now define
$$\hat{h} = \frac{h}{CK^{2k+2}} - K.$$
Then for $\hat{h} \ge 0$ we have
$$\pddt{\hat{h}} \le \triangle \hat{h} - \hat{h}^2.$$
Hence $\hat{h} \le CK$ on $P_{r_k}(X,T)$ from which it follows that $h 
\le CK^{2k+3}$ and 
$$|D^k F| \le CK^{k/2 +1}.$$
\end{proof}

\pagebreak[4]
\bigskip
\setcounter{section}{3}
\setcounter{lemma}{0}
\setcounter{theorem}{0}
\noindent
{\bf 3. Formation of singularities}
\bigskip

In this section, we will prove the main theorem.
Suppose that the flow has a singularity at $X \in M$ at time $T$. We 
suppose that the singularity is rapidly forming, that is, we have the dilation 
invariant estimate
$$|F| \le \frac{C}{T-t},$$
for some constant $C$.  Then Theorem \ref{theoremestimates} gives the dilation 
invariant estimates
$$|D^k F| \le \frac{C_k}{(T-t)^{k/2+1}}.$$
We will perform a blow-up procedure around the point $X$ at time $T$.  Choose a 
small geodesic ball $B(r)$ of radius $r$ around $X$ over which $E$ is trivial.  We
identify $B(r)$ with the ball of radius $r$ in $\mathbf{R}^n$.  In these
coordinates, the connections
$A=A(x,t)$ are given by  matrix-valued one-forms $A = A_p dx^p$. Now choose a
sequence of real numbers
$\lambda_i$ tending to zero, and define a sequence of flows of connection one
forms $A(i) = A_p (i)dx^p$ on balls of increasing size in $\mathbf{R}^n$ by
setting
$$ A_p (i) (x,t) = \lambda_i A_p (\lambda_i x, T+ \lambda_i^2 t)$$
for $x \in B_{1/ \lambda_i}(0)$ and $t \in [- T/ \lambda_i^2 ,0)$.  

We need the following monotonicity formula of Hamilton's \cite{ha} for the 
Yang-Mills flow.
Similar results have been obtained by Naito \cite{na} and Chen and Shen \cite
{cs}.  

\addtocounter{lemma}{1}
\addtocounter{corollary}{1}
\begin{theorem} \label{monotonicitytheorem}
If $A(t)$ solves the Yang-Mills heat flow on $0 \le t < T$ and if $k$ is any 
positive backward solution to the scalar heat equation on $M$ with $\int_M k =1
$, then the quantity
$$Z(t) = (T-t)^2 \int_M |F|^2 k$$
is monotone decreasing in $t$ when $M$ is Ricci parallel (that is, $D_i R_{jk}=0
$) with weakly positive sectional curvatures; while on a general $M$ we have
$$Z(t) \le C Z(\tau) + C(t-\tau)^2 E_0$$
whenever $T-1 \le \tau \le t \le T$, where $E_0$ is the initial energy and $C$ 
is a constant depending only on $M$.  Moreover, if 
$$W(t) = (T-t) \int_M | D_i F^{\al}_{ij \be} + \frac{D_i k}{k}. F^{\al}_{ij 
\be} |^2 k$$
then 
$$\int_{\tau}^T W(t) dt \le C Z(\tau) + C E_0.$$
\end{theorem}

The last statement of Theorem \ref{monotonicitytheorem} gives 
$$ \int_{T-1}^T (T-t)^2 \int_M | \textrm{div} F_A + \frac{D k}{k}. F_A |^2 k \, 
dV \, dt \le C < \infty.$$
Hence for $\epsilon > 0$ there exists $\delta >0$ such that 
$$ \int_{T - \delta}^T (T-t)^2 \int_M |  \textrm{div} F_A + \frac{D k}{k}. F_A 
|^2 k \, dV \, dt \le \epsilon.$$
Now let $k$ be the solution to the backwards heat equation which
becomes a delta function at the point $X$ at time $T$.  In the
geodesic ball, we get
$$ \int_{T - \delta}^T (T-t)^2 \int_{B(r)} |  \textrm{div} F_A + \frac{D k}{k}. 
F_A |^2 k \, dV \, dt \le \epsilon.$$
Let $A(i)$ be the dilated connections and let $k(i)$ be the
dilations of $k$.  Making a change of variables,
$$ \int_{- \delta/ \lambda_i^2}^0 |t|^2 \int_{B_{r/\lambda_i}(0)} | 
\textrm{div}_
{A(i)} F_{A(i)} + \frac{D k(i)}{k(i)} . F_{A(i)}|^2 k(i) \, dV 
\, dt \le \epsilon.$$
Fix a compact time interval $[t' , t'']$ in $(-\infty, 0)$ and a compact 
set $D$ in $\mathbf{R}^n$.  Suppose $i$ is large enough so that $t'$ is 
contained in the interval $( -\delta / \lambda_i^2 , 0)$ and that $D$ is 
contained in the ball $B_{r / \lambda_i}(0)$.  Fix a point $x_0$ in $D$.  
We will now use a theorem of
Uhlenbeck's \cite{u82}, which states that if the $L^n$ norm of the curvature 
is small enough, then there exists a local Coulomb gauge, in which the 
$L^p_1$ norm of the connection matrix is bounded by the $L^p$ norm of the 
curvature, for any $p \ge 2n$.  Since by assumption the curvature is 
uniformly bounded by $C/(T-t)$, we can apply Uhlenbeck's theorem in a 
small enough ball $B$ containing $x_0$.  Thus there exists, for each $i$, 
a gauge 
transformation $s(i)$ over $B$ such that $s(i) \cdot A(i)$ at time $t'$ is 
bounded in $L^p_1$ over $B$ for any $p$.  The bounds do not depend on $i$.  
We can also get uniform bounds on the $C^k$ norms of $s(i) \cdot A(i)$, 
for large enough $i$, at 
$t'$ on a slightly smaller ball $B'$ from the bounds on the derivatives of 
the curvature given in Theorem \ref{theoremestimates}, using 
the argument of Lemma 2.3.11 in \cite{dk}.  Since $s(i) \cdot A(i)$ is a 
solution of the Yang-Mills flow,  we can bound the 
connection matrices and their derivatives uniformly in time on the 
interval $[t', t'']$ over 
the ball $B'$.  By passing to a subsequence, we can get uniform 
convergence 
in $C^{\infty}$.
Since this works for a small ball around any point, we can 
apply the patching argument of \cite{dk}, Corollary 4.4.8, so that after 
passing to another subsequence, $s(i) \cdot A(i)$ converges in 
$C^{\infty}$ on 
$D \times [t',t'']$ to a solution $\overline{A}$ of the Yang-Mills flow.  
Now consider a sequence of increasing compact sets $D_j$ and time 
intervals $[t'_j, t''_j]$ which exhaust $\mathbf{R}^n \times (-\infty,0)$.  
By the usual diagonal argument, after passing to a subsequence, we get a 
sequence of gauge transformations $s(i)$ such that $s(i) \cdot A(i)$ 
converges 
in $C^{\infty}$ on compact sets to a solution $\overline{A}$ of the 
Yang-Mills flow on $\mathbf{R}^n \times (-\infty,0)$.  Also, a subsequence of the
$k(i)$ will converge to $\overline{k}$ on $\mathbf{R}^n \times (-\infty,0)$ given
by
$$\overline{k}(x,t) = \frac{1}{( -4\pi t)^{n/2}} e^{|x|^2/4t}.$$
Letting $\epsilon 
\rightarrow 0$, we get  
$$\textrm{div}_{\overline{A}} F_{\overline{A}} + \frac{x}{2t} . 
F_{\overline{A}} = 0.$$
Hence
$$ \pddt{} \overline{A}^{\al}_{i \be} + \frac{x^p}{2t}(\partial_p
\overline{A}^{\al}_{i \be} - \partial_i \overline{A}^{\al}_{p \be} +
\overline{A}^{\al}_{p \ga} \overline{A}^{\ga}_{i
\be} - \overline{A}^{\al}_{i \ga} \overline{A}^{\ga}_{p
\be}) = 0.$$
Now choose an exponential gauge for $\overline{A}$, in which $x^p \overline{A}^{\al}_{p \be}
= 0$, at time $t=-1$.  Then notice that, from the equation above, the condition
$x^p \overline{A}^{\al}_{p \be} =0$ holds for all $t$.  Then $\overline{A}$
satisfies
$$\pddt{} \overline{A}^{\al}_{i \be} + \frac{x^p}{2t} \partial_p \overline{A}^{\al}_{i \be}
+\frac{1}{2t} \overline{A}^{\al}_{i \be}=0,$$
and by the following lemma $\overline{A}$ is a 
homothetically shrinking soliton.

\addtocounter{theorem}{1}
\addtocounter{corollary}{1}
\begin{lemma} \label{lemmashrinking}
If a flow of one-forms $A_i=A_i(x,t)$ satisfies
$$\pddt{} A_i + \frac{x^p}{2t} \partial_p A_i
+\frac{1}{2t} A_i =0,$$
on $\mathbf{R}^n \times (-\infty,0)$ then 
$$A_i (x,t) = \lambda A_i (\lambda x, \lambda^2 t).$$
\end{lemma}

\begin{proof}
First notice that if $A$ satisfies the dilation invariant condition then by 
differentiating with respect to $\lambda$ and setting $\lambda$ equal to 1, it is 
immediate that it satisfies the differential equation
$$\pddt{} A_i + \frac{x^p}{2t} \partial_p A_i
+\frac{1}{2t} A_i =0,$$
on $\mathbf{R}^n \times (-\infty,0)$.  The lemma then follows from Holmgren's uniqueness
theorem for linear partial differential equations (see for example \cite{ho}, p.125), since
the hyperplanes $\{t=\textrm{constant}\}$ are non-characteristic.
\end{proof}

We will now show that the soliton has non-zero curvature.  First we have 
the 
following `$\epsilon$-regularity' result which is almost the same as one given 
in \cite{gh} for the harmonic map heat flow (a similar result for the Yang-Mills flow is proved in \cite{cs}).  
Write $k_{(X,T)}$ for the backwards solution of the heat equation on $M$ which
becomes a delta function at $X$ at time $T$.

\addtocounter{lemma}{1}
\addtocounter{corollary}{1}
\begin{theorem}
There exist constants $\epsilon >0$ and $\beta >0$ depending only on M and $E_0
$ such that for any $X \in M$, any $T>0$ and any $\al \ge 0$ with $T- \beta \le 
\al <T$ we can find $\rho >0$ and $B < \infty$ such that if $A$ is any solution 
to the Yang-Mills heat flow on $0 \le t \le T$ with energy bounded by $E_0$ and
$$(T-\al)^2 \int_M |F_A (x, \al)|^2 k_{(X,T)} (x, \al) dV \le \epsilon,$$
then $|F_A(x,t)| \le B$ for all $(x,t) \in P_{\rho} (X,T)$.
\end{theorem}
\begin{proof}
We omit the proof, since it is almost identical to that of Theorem 3.2 in \cite
{gh}.
\end{proof}
 
We also have the following corollary.

\begin{corollary}
There exist $\epsilon >0$ and $\beta >0$ depending only on the bundle and $E_0$ 
such that if $A$ is any solution to the Yang-Mills flow on $0 \le t <T$ and if
$$(T-t)^2 \int_M |F_A(x,t)|^2 k_{(X,T)} (x,t) dV \le \epsilon$$
for some $t$ in $T- \be \le t <T$ then $A(x,t)$ extends smoothly on $0 \le t 
\le T$ in some neighbourhood of $X$.
\end{corollary}
\begin{proof}
Apply the above theorem to $\tilde{A}(x,t) = A(x,t - \zeta)$ for $\zeta >0$ and 
let $\zeta \rightarrow 0$. We then get uniform estimates on $|F|$ in some small 
neighbourhood of $X$ for $t<T$.  From Theorem \ref{theoremestimates} we get 
bounds on the higher derivatives, and the result follows.
\end{proof} 

We can now show that the limiting flow $\overline{A}$ does not have zero 
curvature.
Since we have assumed that $X$ is a singular point, we know that for all 
$t$ in $T-\be \le t <T$ we have
$$ (T-t)^2 \int_M |F_{A}(x,t)|^2 k_{(X,T)}(x,t) dV \ge \epsilon $$
for a fixed $\epsilon >0$.  This estimate is dilation invariant.  Now we have 
the bound
$$(T-t)^2 |F|^2 \le C^2.$$
Note that there exists $\rho >0$ such that, with $\overline{k}$ as above,
$$\int_{|x|>\rho \sqrt{|t|}} \overline{k}(x,t) dV \le 
\frac{\epsilon}{2C^2}.$$
Then after making the change of coordinates as before, and taking the 
limit, we see that
$$ |t|^2 \int_{|x| \le \rho \sqrt{|t|}} |F_{\overline{A}} (x,t)|^2 
\overline{k}(x,t) dV \ge \epsilon /2,$$
and hence the soliton has non-zero curvature.

\bigskip
\noindent
\setcounter{section}{4}
{\bf 4. Examples of homothetically shrinking solitons}
\bigskip

We give examples of homothetically shrinking solitons on $\mathbf{R}^n 
\times SO
(n)$ for $5 \le n \le 9$.  It was shown by Gastel \cite{ga} that such 
solitons 
exist in these dimensions.  We follow \cite{gr} and consider 
$SO(n)$-equivariant 
connections given by
$$A_i (x) = -\frac{h(r)}{r^2} \sigma_i (x),$$
where $r= |x|$, $h$ is a real-valued function on $[0, \infty)$ and $\{ \sigma_i 
\}_{i=1}^n$ in $so(n)$ are given by
$$(\sigma_i)^{\al}_{\ \be} = \delta_i^{\al} x^{\be} - \delta_i^{\be} x^{\al}, 
\qquad \textrm{for } 1 \le \al, \be \le n.$$
A long but straightforward calculation shows that the Yang-Mills heat equation 
becomes
$$h_t = h_{rr} + (n-3) \frac{h_r}{r} - (n-2) \frac{h(h-1)(h-2)}{r^2}.$$
Define a function 
$$\phi_(\rho) = \frac{\rho^2}{ a_n \rho^2 + b_n},$$
where $a_n$ and $b_n$ are positive constants given by
$$a_n = \frac{\sqrt{n-2}}{2\sqrt{2}} \qquad \textrm{and} \qquad b_n = \frac{1}
{2} (6n -12 - (n+2)\sqrt{2n -4}).$$
Notice that $b_n$ is positive if and only if $5 \le n \le 9$. 
Now for $t$ in $(-\infty ,0)$, let
$$h(r,t) = \phi(\rho), \qquad \textrm{where} \qquad \rho = \frac{r}{\sqrt{-
t}}.$$
Then it is easy to check that this $h(r,t)$ gives a smooth solution
$A=A(x,t)$ of the  Yang-Mills heat equation on $\mathbf{R}^n \times
(-\infty, 0)$ and that
$$A_i(x,t) = \lambda A_i(\lambda x , \lambda^2 t).$$

\end{document}